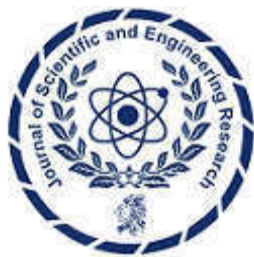

# On a generalized theorem of de Bruijn and Erdös in d-dimensional Fuzzy Linear Spaces

Hasan Keleş

Karadeniz Technical University, Faculty of Science, Department of Mathematics, Ortahisar, Trabzon, Turkey

**Abstract** In this study we follow a new framework for the theory that offers us, other than traditional, a new angle to observe and investigate some relations between finite sets, $F$-lattice $L$ and their elements.
The theory is based on the *Fuzzy Linear Spaces* (*FLS*) $S = (N, D)$. In this case, to operate on these spaces the necessary preliminaries, concepts and operations in lattices relative to *FLS* are introduced. Some definitions, such that *k-fuzzy point*, *k-fuzzy line* are given. Then we correspond these definitions to the definitions in usually linear spaces. We investigate some combinatorics properties of *FLS*. In some examples in the case where $|L| = 3$ .
We see some differences. In general, taking an ordered lattice $L_n = \{0, a_1, a_2, \ldots, a_n, 1\}$ we observe how some combinatorics formulas and properties are changed. In *FLS* the dimension concept is a set. We produce some general formulas by using some trivial examples. Furthermore, we generalize de Bruijn-Erdös Theorem in [4].

**Introduction**

k-point, k-line for Linear Spaces, d-dimensional Linear Spaces were studied by some authors like Batten [5] and Barwick [6]. Here, we give a very short proof to well-known the theorem of de Bruijn and Erdös [4,5][↑]. And also, we have been collected all them from the above papers and from [1,2,4,5].
In this paper, we extended the Theorem de Bruijn and Erdös. For this we have to give.

**Definition 1.** Let $S = (N, D)$ be a *FLS* and $X \subset N$. The set

$$\left\{ x \in N : \forall x_1, \ldots, x_n \in X, \exists d \in D, \bigwedge_{\substack{i=1 \\ k, m \in \mathbb{N} \text{ ve } 2 \leq k \leq m}}^{k} d(x_i) \wedge d(x) \neq \theta \right\}$$

is called *closure* of $X$ and denoted by $\langle X \rangle$.

In any $S = (N, D)$ *FLS*, $\langle \emptyset \rangle = \emptyset$, $\langle \{x\} \rangle = \{x\}$ and $\langle S \rangle = S$.

---

• $|L|$: Number elements of $L$.
↑ De Bruijn and Erdös (1948). Sometimes called the de Bruijn-Erdös and Hanani theorem because of Hanani (1955).

*Journal of Scientific and Engineering Research*

If $\langle X \rangle = B$ then we say that $X$ generates $B$.

**Example 2.** Let $S = (N, D)$ FLS, where $N = \{x, y, z\}$ and $D = \{d_1 = \{1,1,0\}, d_2 = \{0,1,1\}, d_3 = \{1,0,1\}\}$. For $X = \{x, y\}$, $\langle X \rangle = d_1$, which is only one line.

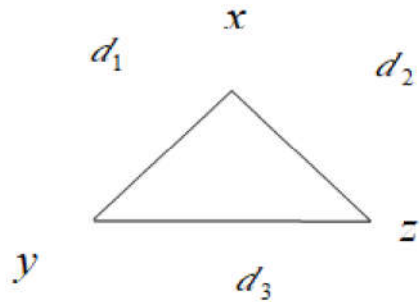

If $L_1 = \{0, a_1, 1\}$, $\langle X \rangle = d_1 \cup \{1, a_1, 0\} \cup \{a_1, 1, 0\} \cup \{a_1, a_1, 0\}$ is not one line.

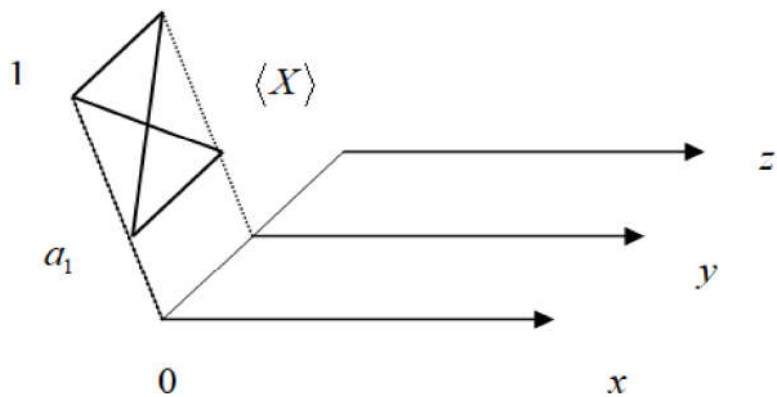

**Definition 3.** Let $S = (N, D)$ be a FLS. Then any point $x \in N$ is called *k-fuzzy point* if $\bigwedge_{i=1, d_i \in D}^{k} d_i(x) \neq \theta$.

$\{x\}$ is *0-fuzzy point* for $L = \{0,1\}$. But $\{x\}$ is *1-fuzzy point* for $L_n, n \geq 2$.

**Definition 4.** Let $S = (N, D)$ be a FLS. Then a line $d \in D$ is called *k-fuzzy line* if $\bigwedge_{i=1, x_i \in N}^{k} d(x_i) \neq \theta$.

**Lemma 5.** Let $S = (N, D)$ be a FLS and any line $d \in D$ be a k-fuzzy line. Then the number of k-fuzzy line is $(n+1)^k$.

**Proof.** There are k points $x_1, \ldots, x_k$ on each k-fuzzy line and $d(x_j) = t$ where $t = a_1, \ldots, a_n, 1$ and $j = 1, \ldots, k$. Then the number of k-fuzzy line is $(n+1)^k$.

**Lemma 6.** Let $S = (N, D)$ be a FLS and any point $x \in N$ be a k-fuzzy point. Then the number of k-fuzzy points is $\prod_{j=1}^{k} (n+1)^{v_j}$ where $v_j = |\{x | d_j(x) \neq \theta, x \in N\}|$.

**Proof.** If there are $v_j$ points on each line $d_j$ from Lemma 5 the number of such line $d_j$ just $(n+1)^{v_j}$, where $j = 1, \ldots, k$. And furthermore since $x$ is a k-fuzzy point then the total number of k-fuzzy points is $\prod_{j=1}^{k} (n+1)^{v_j}$.

**Theorem** (de Bruijn-Erdös) [5]. Let *S* be any finite linear space with $b = |D| > 1, |N| = v$. Then

i. $b \geq v$,

ii. If $b = v$, any two lines have point in common. In case (2) either one line has $v - 1$ points and all others have two points, or every line has $k + 1$ points and every point is on $k + 1$ lines, $k \geq 2$.

If any point of *S* has k-fuzzy point then the following proposition will give:

**Proposition 7.** Let $S = (N, D)$ be a *FLS* such that $|S| = m$, any point $x \in N$ is k-fuzzy point, and $v_j = |\{x | d_j(x) \neq \theta, x \in N\}|$. Then

$$|D| = \frac{\log_{(n+1)} m}{v_j}, n \geq 1.$$

**Proof.** $|S| = m$ by [2].

$$= \prod_{j=1}^{|D|} (n+1)^{v_j}$$

$$= \underbrace{(n+1)^{v_j} \ldots (n+1)^{v_j}}_{|D|\text{-time}}$$

$$= (n+1)^{|D|v_j}$$

$\log_{(n+1)} m = |D| v_j \log_{(n+1)} (n+1)$

$$|D| = \frac{\log_{(n+1)} m}{v_j}, n \geq 1.$$

We now extend the Theorem of de Bruijn-Erdös:

**Theorem** (Hasan KELEŞ). Let $S = (N, D)$ be any finite *FLS* such that with $b = |D| > 1, |N| = v \geq 3$. Then $b \geq v$ and any two lines have a point in common. Furthermore, either just one of the lines in *D* is a $(v - 1)$-fuzzy line and others are *2*-fuzzy lines, or every line is a $(k+1)$-fuzzy line and every point is $\left[ \prod_{j=1}^{k+1} (|L| - 1)^{v_j} \right]$ -fuzzy point, $k \geq 2$.

**Proof.** The inequality $b \geq v$ is obvious. The case where $|L| = 2$ it is the theorem *de Bruijn Erdös'*. It is clear that $b > v$. The fact that any two lines have a point in common is obtained from the definition of *FLS*. If one of the lines in *D* is $(v-1)$-fuzzy line then $\bigwedge_{i=1}^{v-1} d(x_i) \neq \theta$ and $d(x_v) = \theta$. So $d(x_v) \bigwedge_{i \in \{1,\ldots,v-1\}} d'(x_i) \neq \theta$ for $\forall d' \neq d, d' \in D$ from the definition of *FLS*. Therefore lines $d'$ are 2-fuzzy lines.

If one of the lines in *D* is not a $(v-1)$-fuzzy line then the other are not 2-fuzzy lines. Therefore all of them are $(k+1)$-fuzzy lines where $k \geq 2$. So $\bigwedge_{i=1}^{k+1} d(x_i) \neq \theta$. Line *d* has points $(k+1)$. Therefore any point $x_i \in N$ is a $\left[ \prod_{j=1}^{k+1} (|L|-1)^{v_j} \right]$ -fuzzy point from Lemma 6.

*Journal of Scientific and Engineering Research*

**References**
[1]. Keles, H.,(1996). *d*-dimensional Linear Spaces, University of Sakarya, Turkey, 47p.
[2]. Keles, H., (2017). On General Construction of d-dimensional Linear Spaces, Journal of Scientific and Engineering Research, vol.4 pp.365.
[3]. Koen Thas. Automorphisms and combinatorics of finite generalized quadrancles. 440p.
[4]. L. M. Batten, A. Beutelspacher., (1993). The theory of finite linear spaces. Combinatorics of Points and Lines. Cambridge University Press, 214p.
[5]. L. Batten., (1986). Combinatorics of finite geometries. Cambridge University Press, 172p.
[6]. S. Barwick., (1994).Substructures of finite geometries. The University of London, 113p.
[7]. Keleş, H., (2017). On some properties combinatorics of Graphs in the d-dimensional FLS, Journal of Scientific and Engineering Research, vol.4 pp.93,
[8]. Keleş, H., (2005). On the d-dimensional Fuzzy Linear Spaces, SamTa05, Sampling theory and applications International Workshop, 45.
[9]. Keleş, H., (2004), d-dimensional Fuzzy Linear Spaces, International workshop on global analysis, 20.
[10]. Keleş, H., (2006), On Some Numbers Related to the Differential Equation System of d-Dimensional Fuzzy Linear Spaces (FLS), Mathematical Methods in Engineering, International Symposium, 33.
[11]. Diestel, R., (2006), Graph Theory", Springer, 3$^{rd}$ Edition, 1-33, Hamburg, Germany.
[12]. Keleş, H., (2016), On Measurement Assessment and Division, Journal of Scientific and Engineering Research, 3(6): 233-237.